\documentclass[10pt,twoside,english]{myamsart}
\advance\oddsidemargin by -1.0cm
\advance\evensidemargin by -1.0cm
\advance\topmargin by -1.0cm 
\usepackage{amssymb}
\usepackage{babel}
\usepackage{amstext}
\usepackage{amscd}   
\usepackage{epsfig}  
\usepackage{rotating}
\let\mathg\mathfrak
\theoremstyle{plain}

\newtheorem{lem}{Lemma}[section]
\newtheorem{thm}{Theorem}[section]
\newtheorem{prop}{Proposition}[section]
\theoremstyle{definition}
\newtheorem{exa}{Example}[section]
\newtheorem{NB}{Remark}[section]

\renewcommand{\labelenumi}{$(\theenumi)$}
\newcommand{\bdm}{\begin{displaymath}}
\newcommand{\edm}{\end{displaymath}}
\newcommand{\be}{\begin{equation}}
\newcommand{\ee}{\end{equation}}
\newcommand{\ba}[1]{\begin{array}{#1}}
\newcommand{\ea}{\end{array}}
\newcommand{\bea}[1][]{\begin{eqnarray#1}}
\newcommand{\eea}[1][]{\end{eqnarray#1}}
\newcommand{\btab}{\begin{tabular}}
\newcommand{\etab}{\end{tabular}}

\newcommand{\Id}{\ensuremath{\mathrm{Id}}}

\newcommand{\R}{\ensuremath{\mathbb{R}}}

\newcommand{\Z}{\ensuremath{\mathbb{Z}}}

\newcommand{\Ric}{\ensuremath{\mathrm{Ric}}}
\newcommand{\Scal}{\ensuremath{\mathrm{Scal}}}

\newcommand{\diag}{\ensuremath{\mathrm{diag}}}

\newcommand{\SU}{\ensuremath{\mathrm{SU}}}
\newcommand{\U}{\ensuremath{\mathrm{U}}}

\newcommand{\SO}{\ensuremath{\mathrm{SO}}}
\newcommand{\Spin}{\ensuremath{\mathrm{Spin}}}

\begin{document}
\def\haken{\mathbin{\hbox to 6pt{%
                 \vrule height0.4pt width5pt depth0pt
                 \kern-.4pt
                 \vrule height6pt width0.4pt depth0pt\hss}}}
    \let \hook\intprod
\setcounter{equation}{0}
%
%
\thispagestyle{empty}
%
\date{\today}
\title{Almost Contact Manifolds, Connections with Torsion, and Parallel Spinors}
%
%
%
\author{Thomas Friedrich and Stefan Ivanov}
\address{\hspace{-5mm} 
Thomas Friedrich\newline
Institut f\"ur Reine Mathematik \newline
Humboldt-Universit\"at zu Berlin\newline
Sitz: WBC Adlershof\newline
D-10099 Berlin, Germany\newline
{\normalfont\ttfamily friedric@mathematik.hu-berlin.de}}
\address{\hspace{-5mm} 
Stefan Ivanov\newline
Faculty of Mathematics and Informatics\newline
University of Sofia ``St. Kl. Ohridski''\newline
blvd. James Bourchier 5\newline
1164 Sofia, Bulgaria\newline
{\normalfont\ttfamily ivanovsp@fmi.uni-sofia.bg}}
\thanks{Supported by the SFB 288 "Differential geometry
and quantum physics" of the DFG and the European Human Potential Program EDGE, 
Research Training Network HPRN-CT-2000-00101. S.Ivanov thanks ICTP for the
support and excellent environments.}
\keywords{Contact structures, string equations}  
\begin{abstract}
We classify locally homogeneous quasi-Sasakian manifolds in dimension five 
that admit a parallel spinor $\psi$ of algebraic type $F \cdot \psi = 0$ with 
respect to the unique connection $\nabla$ preserving the quasi-Sasakian 
structure and with totally  skew-symmetric torsion. We introduce a certain
conformal transformation of almost contact metric manifolds and discuss a
link between them and the dilation function in $5$-dimensional string theory.
We find natural conditions implying conformal invariances of parallel
spinors. We present topological obstructions to the existence of parallel
spinors in the compact case.
\end{abstract}
\maketitle
\pagestyle{headings}
%
%
%
\section{Introduction}\noindent
The basic model in type II string theory is a $6$-tuple 
$(M^n,g,\nabla,T,\Phi, \psi)$ consisting of a Riemannian metric $g$,
a metric connection $\nabla$ with totally skew-symmetric torsion form $T$, 
a dilation function $\Phi$ and a spinor field $\psi$. The string 
equations can then be written in the following form (see \cite{Stro} and
\cite{IP, Friedrich&I}):
\bdm
\mbox{Ric}^{\nabla} + \frac{1}{2} \cdot \delta^g(T) + 2 \cdot \nabla^g(d \Phi)
\ = \ 0, \quad \delta^g(T) \ = \ 2 \cdot (\mbox{grad}(\Phi) \haken T),
\edm
\bdm
\nabla \psi \ = \ 0, \quad (2 \cdot d \Phi - T) \cdot \psi \ = \ 0 \, .
\edm
If the dilation function is constant the string equations simplify:
\bdm
\mbox{Ric}^{\nabla} \ = \ 0, \quad  \delta^g(T) \ = \ 0, 
\quad \nabla \psi \ = \ 0, \quad T \cdot \psi \ = \ 0 \, .
\edm
In fact, the bosonic part is taken with one-loop contribution. With two-loop
contribution it receives additional terms involving quadratic terms of the 
curvature (see \cite{HP}). The fermionic part, i.e., the Killing spinor 
equation, is responsible for the preserved supersymmetries. The number of
preserved supersymmetries depends on the number of parallel spinors. In this
paper we concentrate our attention mainly on the geometry of the solutions 
of the Killing spinor equations in dimension five. In this dimension 
the $\nabla$-parallel spinor field $\psi$ defines, via the formulas
\bdm
\xi \cdot \psi \ = \ i \cdot \psi, \quad - \, 2 \cdot \varphi(X) \cdot \psi + 
\xi \cdot X \cdot \psi \ = \ i \cdot X \cdot \psi \, ,
\edm
an almost contact metric structure $(M^5, g , \xi , \eta , \varphi)$, which
is preserved by the connection $\nabla$. A simple algebraic computation 
yields the equations
\bdm
\varphi(\xi) \ = \ 0, \quad \varphi^2 \ = \ - \, \mbox{Id} + \eta \otimes \xi,
\quad g(\varphi(X),\varphi(Y)) \  = \ g(X,Y) - \eta(X) \cdot \eta(Y) \, .
\edm
However, the contact condition $\eta \wedge (d \eta)^2 \neq 0$ need not be 
satisfied in general. This suggests that in dimension five the solutions of 
the type II string equations are related to suitable almost contact 
metric structures. 
In fact, not any such geometric structure is admissible. The first condition 
is that there should exist a connection with totally skew-symmetric torsion 
and preserving the almost contact metric structure. In the paper 
\cite{Friedrich&I} we proved that an almost contact metric 
structure admits such a connection 
if and only if the Nijenhuis tensor $N$ is totally skew-symmetric and if the 
vector field $\xi$ is a Killing vector field. In this case the connection is 
unique and we derived a formula for the torsion form $T$. Suppose now that 
we fix an almost contact metric structure of that type. Then we have a 
reduction of the frame bundle  $\mathcal{F}(M^5)$ to the subgroup 
$\U(2) \subset \SO(5)$. 
But the isotropy group of a spinor in the $5$-dimensional spin representation 
$\Delta_5$ does not coincide with $\U(2)$. Therefore the existence of a 
$\nabla$-parallel spinor field imposes a second condition on the almost 
contact metric structure under consideration. The almost contact metric 
structure splits the spinor bundle $\Sigma$ of $M^5$ into two $1$-dimensional
bundles $\Sigma^{\pm}$ and into one $2$-dimensional bundle $\Sigma^2$. Since
the unique connection $\nabla$ preserves the almost contact metric structure,
it preserves the decomposition of the spinor bundle, too. Consequently, we 
obtain two different integrability conditions for the $\nabla$-parallel 
spinor. The case that the spinor field is a section in one of the 
$1$-dimensional subbundles $\Sigma^{\pm}$ means that the connection 
$\nabla$ has a reduction to the subgroup $\SU(2)$. Equivalently, the Ricci 
tensor of the connection has to satisfy certain algebraic relations. This 
case was studied in \cite{Friedrich&I}. In
particular, there are compact Sasakian manifolds with $\nabla$-parallel spinor
fields of that algebraic type. The aim of this paper is to study the 
almost contact metric structures in dimension five admitting a 
$\nabla$-parallel spinor field inside the $2$-dimensional bundle $\Sigma^2$.  
We classify all locally homogeneous quasi-Sasakian structures with a 
$\nabla$-parallel spinor field in the bundle $\Sigma^2$. Furthermore, we discuss the second Killing equation for spinors in the bundle $\Sigma^2$ 
involving a dilation function and show that it has to be constant. Concerning 
solutions of both Killing spinor equations in the bundle $\Sigma^{\pm}$, 
we introduce certain transformations depending on a real function (special conformal transformations) of an almost contact metric structure and show that any
solution to both Killing equations in the bundle $\Sigma^{\pm}$ is invariant under these transformations in dimension five. We discuss the close relationship 
between the dilation function and the Lee form of the structure and show that 
the dilation function can be interpreted as a conformal factor. In the regular case, we find that the second Killing equation implies that any parallel spinor
is projectable. This allows us to list all compact regular solutions to both equations in dimension five. Finally we prove, in any odd dimension, a generalization of Tachibana's Theorem for 
harmonic $1$-forms on compact quasi-Sasakian manifolds in the presence of 
some special $\nabla$-parallel spinor.

\section{Contact connections with parallel spinors}\noindent
%
We start with some basic definitions in contact geometry, on which
the book \cite{Blair} or the article \cite{CG} may serve as a general
reference.
An \emph{almost contact metric structure} consists of an odd-dimensional 
manifold $M^{2k+1}$ equipped with a Riemannian
metric $g$, a vector field $\xi$ of length one, its dual $1$-form $\eta$ 
as well as  an endomorphism $\varphi$ of the tangent bundle 
such that the algebraic relations
\bdm
\varphi(\xi) \ = \ 0, \quad \varphi^2 \ = \ - \Id + \eta \otimes \xi, \quad
g(\varphi(X),\varphi(Y)) \  = \ g(X,Y) - \eta(X) \cdot \eta(Y) 
\edm
are satisfied. If $\eta \wedge (d \eta)^k \neq 0$, we have a \emph{contact
manifold}. The \emph{fundamental form} $F$ and the \emph{Nijenhuis tensor} $N$
of an almost contact metric structure are defined by the formulas
\bdm
F(X,Y) \ := \ g(X, \varphi(Y))\, , \quad 
N(X,Y) \ := \ [\varphi, \varphi](X,Y) + d \eta (X,Y) \cdot \xi \, .
\edm
It is clear that $\eta \wedge F^k \neq 0$. There are many special types of 
almost contact metric structures. We introduce
those appearing in this paper. An almost
contact metric structure is called \emph{normal} if its Nijenhuis tensor
vanishes, $N = 0$. A \emph{quasi-Sasakian} structure is a normal
almost contact metric structure with closed fundamental form, 
$N = 0$ and $dF = 0$. The vector field $\xi$ of a quasi-Sasakian structure is 
automatically a Killing vector field. In fact, in Section $3$ of the present
paper we will prove a more general result. A normal, almost contact metric 
structure with the property that the derivative $d\eta$ of the contact form 
is proportional to the fundamental form are called \emph{$\alpha$-Sasakian},  
$N = 0$, $d \eta = \alpha \cdot F$ and $\alpha$ constant. 
Any $\alpha$-Sasakian structure
is quasi-Sasakian. Finally, \emph{Sasakian}
manifolds are characterized by the following two integrability conditions
\bdm
N \ = \ 0 \, , \quad d \eta \ = \ 2 \cdot F \, .
\edm

\noindent
Now we study the spin geometry of the $5$-dimensional
local model $\R^5$ in detail. Let us fix an orthonormal basis 
$e_1,\,\ldots, e_5$ and consider the almost contact metric structure given by 
the vector $\xi := e_5$ and the skew-symmetric endomorphism
$\varphi := - e_1 \wedge e_2 - e_3 \wedge e_4$. The fundamental form $F$ is 
\bdm
F \ = \ e_1 \wedge e_2 + e_3 \wedge e_4 \, . 
\edm
The subgroup of $\SO(5)$ preserving the $(g, \xi, \eta, \varphi)$-structure is 
isomorphic to the group $\U(2)$. A form $\sum x_{ij} \cdot e_i \wedge e_j$ 
belongs to the Lie algebra $\mathg{u}(2)$ if and only if
\bdm
x_{14} + x_{23} \ = \ 0, \quad x_{13} - x_{24} \ = \ 0, \quad x_{i5} \ = \ 
0 \, . 
\edm
Denote by $\Delta_5$ the $5$-dimensional spin representation. Vectors act
on $\Delta_5$ by Clifford multiplication and we will use the following matrix 
representation (see \cite{Fri2}):
\bdm
e_1 \ = \ \left[\ba{cccc}0 & 0 & 0 & i \\ 0 & 0 & i & 0\\  0 & i & 0 & 0\\ 
i & 0 & 0 & 0  \ea\right], \quad e_2 \ = \ \left[\ba{cccc}0 & 0 & 0 & -1 \\ 0
& 0 & 1 & 0\\  0 & -1 & 0 & 0\\ 1 & 0 & 0 & 0  \ea\right], \quad  e_3 \ = \
\left[\ba{cccc}0 & 0 & -i & 0 \\ 0 & 0 & 0 & i\\  -i & 0 & 0 & 0\\ 0 & i & 0 &
  0  \ea\right] \, ,
\edm
\bdm
 e_4 \ = \ \left[\ba{cccc}0 & 0 & 1 & 0 \\ 0 & 0 & 0 & 1\\  -1 & 0 & 0 & 0\\ 0
   & -1 & 0 & 0  \ea\right], \quad  e_5 \ = \ \left[\ba{cccc}i & 0 & 0 & 0 \\
   0 & i & 0 & 0\\ 0 & 0 & -i & 0\\ 0 & 0 & 0 & -i  \ea\right], \quad
F \ = \ \left[\ba{cccc}2i & 0 & 0 & 0 \\ 0 & -2i & 0 & 0\\  0 & 0 & 0 & 0\\ 
0 & 0 & 0 & 0  \ea\right]\, .
\edm
The fundamental form acts on $\Delta_5$ with eigenvalues $(2i, -2i, 0, 0)$.
Consequently, $\Delta_5$ splits as a $\mathg{u}(2)$-representation into two
$1$-dimensional representations $\Delta_5^{\pm}$ and one $2$-representation
$\Delta_5^2$. These spaces are defined by the conditions
\bdm
\Delta_5^{\pm} \ = \ \big\{ \psi \in \Delta_5 : \ F \cdot \psi \ = \ \pm 2i 
\cdot \psi \big\}, \quad \Delta_5^2 \ = \  \big\{ \psi \in \Delta_5 : \ F
\cdot \psi \ = \ 0 \big\} \, .
\edm 
If we split $\Delta_5 = \Delta_4^+ \oplus \Delta_4^-$ as a
$\mathg{spin}(4)$-representation, then we will obtain $\Delta_4^+ = \Delta_5^+ \oplus 
\Delta_5^-$ and $\Delta_4^- = \Delta_5^2$. The isotropy subgroup of a spinor 
in $\Delta_5^{\pm}$ is isomorphic to the group $\SU(2)$ and its Lie algebra is given by the formulas:
\bdm
x_{14} + x_{23} \ = \ 0, \quad x_{13} - x_{24} \ = \ 0, \quad 
x_{12} + x_{34} \ = \ 0, \quad x_{i5} \ = \ 0 \, . 
\edm
On the other hand, the isotropy subgroup inside $\U(2)$ of a spinor in 
$\Delta_5^2$ is the diagonally embedded subgroup $\U(1)$:
\bdm
x_{12}  \ = \ x_{34}, \quad x_{13} \ = \ x_{14} \ = \ x_{23} \ = \ 
x_{24} \ = x_{i5} \ = \ 0 \, .
\edm
Remark that the diagonally embedded subgroup $\U(1) \subset \U(2)$ induces the
trivial homomorphism on the fundamental groups
\bdm
\pi_1(\U(1)) \longrightarrow \pi_1(\SO(4)) \ = \ \Z_2 
\edm
and, consequently, it lifts into the spin group $\Spin(4)$. An easy algebraic 
computation proves the following
\begin{lem}
Let $\omega \in \mathg{u}(2)$ be a $2$-form in the Lie algebra of the group 
$\U(2)$. Then the following conditions are equivalent:
\begin{enumerate}
\item There exists a nontrivial spinor $\psi \in \Delta_5^2$ such that
  $\ \omega \cdot \psi =  0 $.
\item For any spinor $\psi \in \Delta_5^2$ the Clifford product $\ 
\omega \cdot \psi =  0 $ vanishes.
\item $\omega$ is proportional to the fundamental form, $\ \omega =  
a \cdot F$.
\end{enumerate}
\end{lem}
\noindent
Let $(M^5, g, \xi, \eta, \varphi)$ be an almost contact metric manifold and
denote by $\mathcal{F}(M^5)$ its Riemannian frame bundle. The almost contact
metric structure defines a $\U(2)$-reduction $\mathcal{R} \subset
\mathcal{F}(M^5)$ of the principal $\SO(5)$-bundle $\mathcal{F}(M^5)$. We fix
a spin structure and we denote by $\Sigma$ the spinor bundle. Again, the 
fundamental form of the almost contact metric structure splits the spinor
bundle into two $1$-dimensional subbundles $\Sigma^{\pm}$ and one
$2$-dimensional bundle $\Sigma^2$. Consider a metric connection $\nabla$
preserving the almost contact metric structure, i.e., a connection in
$\mathcal{R}$. Since the connection preserves the
decomposition of the spinor bundle too, we can study the
problem whether or not there are $\nabla$-parallel spinors in any of these
three subbundles of the spinor bundle. In case the $\nabla$-parallel
spinor field is a section in one of the bundles $\Sigma^{\pm}$, the connection
reduces to the subgroup $\SU(2)$ and we obtain only one necessary and 
sufficient condition: 
\bdm
R^{\nabla}_{XY12} + R^{\nabla}_{XY34}  \ = \ 0 \, .
\edm
Moreover, if the connection $\nabla$ has totally skew-symmetric torsion
$T$, then the latter condition becomes equivalent to a certain relation between
the Ricci tensor $\Ric^{\nabla}$ and the exterior derivative $dT$ of the
torsion form (see \cite{Friedrich&I}, Proposition 9.1). In particular, for
Sasakian manifolds and their unique connection with totally skew-symmetric
torsion the integrability condition is
\bdm
\Ric^{\nabla} \ = \ 4 \cdot (g \, - \, \eta \otimes \eta) \, .
\edm 
Compact examples of this type are known (see
\cite{Friedrich&I}, Proposition 9.2 and Remark 9.3).\\ 

\noindent
In this paper we shall study the integrability condition
for $\nabla$-parallel spinors $\psi \in \Gamma(\Sigma^2)$ in the
$2$-dimensional subbundle,
\bdm
\nabla \psi \ = \ 0, \quad F \cdot \psi \ = \ 0 \, .
\edm 
Then the connection $\nabla$ reduces to a $\U(1)$-subbundle $\mathcal{R}_0
\subset \mathcal{R}$. Using an adapted frame in $\mathcal{R}_0$ the abelian 
connection $\nabla$ is given by a $1$-form A:
\bdm
\nabla e_1 \ = \ A \cdot e_2, \quad \nabla e_2 \ = \ - A \cdot e_1,
\quad \nabla e_3 \ = \ A \cdot e_4, \quad \nabla e_4 \ = \ - A \cdot e_4,
\quad  \nabla e_5 \ = \ 0 \, .
\edm
The curvature form $\Omega^A := dA$ of the $\mathcal{R}_0$-connection A is a
well defined $2$-form on $M^5$ and the curvature tensor $R^{\nabla}$
is given by the formulas
\bdm
R^{\nabla}(X,Y)e_1\ = \ \Omega^A(X,Y)e_2, \quad R^{\nabla}(X,Y)e_2\ = \ 
- \, \Omega^A(X,Y)e_1, \quad R^{\nabla}(X,Y)e_3\ = \ \Omega^A(X,Y)e_4, 
\edm
\bdm
R^{\nabla}(X,Y)e_4\ = \ - \, \Omega^A(X,Y)e_3, \quad 
R^{\nabla}(X,Y)e_5\ = \ 0 \, .
\edm
Consequently, we obtain the following condition on the curvature for
the existence of $\nabla$-parallel spinors in the $2$-dimensional bundle
$\Sigma^2$.
\begin{thm}\label{thm-curv} 
Let $(M^5, g, \xi, \eta, \varphi)$ be a simply connected,  almost
  contact metric spin manifold and let $\nabla$ be a connection preserving the
structure. Then the following conditions are equivalent:
\begin{enumerate}
\item There exists a nontrivial, $\nabla$-parallel spinor field in the
  subbundle $\Sigma^2$,
\bdm
\nabla \psi \ = \ 0, \quad F \cdot \psi \ = \ 0 \, .
\edm
\item There are two $\nabla$-parallel spinor fields in the subbundle 
$\Sigma^2$.
\item The curvature tensor $R^{\nabla} : \Lambda^2(M^5) \rightarrow \Lambda^2(M^5)$ is
  given by the formula
\bdm
R^{\nabla}(\alpha) \ = \ (\Omega^A , \alpha) \cdot F \, ,
\edm
where $\Omega^A$ is a closed $2$-form.
\end{enumerate}
\end{thm} 
\begin{NB} A connection $\nabla$ with a parallel spinor in one of the bundles
$\Sigma^{\pm}$ and a parallel spinor in the bundle $\Sigma^2$ is flat. Indeed,
the condition $R^{\nabla}_{XY12} + R^{\nabla}_{XY34} = 0$ yields that the 
curvature form $\Omega^A$ vanishes. If $M^5$ is compact, then it is a 
$5$-dimensional compact Lie group (see \cite{SSTP}).
\end{NB}
\noindent
A direct computation of the $\nabla$-Ricci tensor yields its 
components:
\bdm
R^{\nabla}_{11} \ = \ - \, \Omega^A(e_1,e_2), \quad 
R^{\nabla}_{13} \ = \ - \, \Omega^A(e_1,e_3), \quad 
R^{\nabla}_{14} \ = \  \Omega^A(e_1,e_3),
\edm
\bdm
R^{\nabla}_{22} \ = \ - \, \Omega^A(e_1,e_2), \quad 
R^{\nabla}_{23} \ = \ - \, \Omega^A(e_2,e_4), \quad 
R^{\nabla}_{24} \ = \  \Omega^A(e_2,e_3),
\edm
\bdm
R^{\nabla}_{31} \ = \  \Omega^A(e_2,e_3), \quad 
R^{\nabla}_{32} \ = \ - \, \Omega^A(e_1,e_3), \quad 
R^{\nabla}_{33} \ = \ - \, \Omega^A(e_3,e_4),
\edm
\bdm
R^{\nabla}_{41} \ = \  \Omega^A(e_2,e_4), \quad 
R^{\nabla}_{42} \ = \  -\, \Omega^A(e_1,e_4), \quad 
R^{\nabla}_{44} \ = \  - \, \Omega^A(e_3,e_4),
\edm
\bdm
R^{\nabla}_{51} \ = \  \Omega^A(e_2,e_5), \quad
R^{\nabla}_{52} \ = \  - \, \Omega^A(e_1,e_5), \quad
R^{\nabla}_{53} \ = \  \Omega^A(e_4,e_5), \quad
R^{\nabla}_{54} \ = \  - \, \Omega^A(e_3,e_5)\,. 
\edm
In particular, the two equations $\Ric^{\nabla} = 0 , \nabla \psi =
0$ with a spinor field $\psi \in \Gamma(\Sigma^2)$ in the $2$-dimensional 
bundle have a solution only in case the connection $\nabla$ is flat.
\begin{thm}  
Let $(M^5, g, \xi, \eta, \varphi)$ be an  almost contact metric 
spin manifold and let $\nabla$ be a connection preserving the
structure. If the Ricci tensor of $\nabla$ vanishes
and the connection $\nabla$ admits a parallel spinor such that $F \cdot \psi 
= 0$, then the connection $\nabla$ is flat.
\end{thm}
%
\section{Normal almost contact metric structures with $\nabla$-parallel
spinors}\noindent
%
We suppose that the almost contact metric spin manifold $(M^{2k+1}, g, \xi,
\eta, \varphi)$ admits a connection $\nabla$ with totally skew-symmetric
torsion $T$ and preserving the structure. A connection of this
type exists if and only if the vector field $\xi$ is a Killing vector field
and the Nijenhuis tensor $N$ treated as a $(0,3)$-tensor is totally
skew-symmetric (see \cite{Friedrich&I}). Moreover, the connection is unique
and its torsion form $T$ is given by the formula
\bdm
T \ = \ \eta \wedge d\eta + d^{\varphi}F + N - \eta \wedge (\xi \haken N) \, ,
\edm
where $d^{\varphi}F$ denotes the $\varphi$-twisted exterior differential
of the fundamental form $F$. In particular, we have
\bdm
g(\nabla_X^g \xi, Z) \ = \ \frac{1}{2} T(\xi, X, Z), \quad d \eta \ = \
\xi \haken T, \quad \xi \haken d \eta \ = \ 0 \, , 
\edm
where $\nabla^g$ is the Levi-Civita connection of the Riemannian manifold. If the Nijenhuis tensor $N$ is totally skew-symmetric then the Killing condition
for $\xi$ becomes equivalent to the equation $\xi \haken d F = 0$. In fact, we have
\begin{prop}
Let $(M^{2k+1}, g, \xi, \eta, \varphi)$ be an almost contact metric manifold
with totally skew-symmetric Nijenhuis tensor $N$. The vector field 
$\xi$ is Killing if and only if $\xi \haken d F = 0$.
\end{prop}
\begin{proof} In general, the vector field $\xi$ is a Killing field if and 
only if
\bdm
\nabla_X \xi \ = \ \frac{1}{2} \cdot X \haken d \eta
\edm
holds (see \cite{Blair}, page 64). Using the equation (see \cite{Blair}, page 53)
\bdm
2 \cdot g\big((\nabla_X \varphi)Y, \xi\big) \ = \ - \, dF(X,Y,\xi) + N(Y, \xi, 
\varphi(X)) + d \eta(\varphi(Y),X)
\edm
we immediately obtain the following expression for $X \haken d \eta - 
2 \cdot \nabla_X \xi$:
\bdm
g\big( X \haken d \eta - 2 \cdot \nabla_X \xi, \varphi(Y)\big) \ = \ 
- \, dF(X,Y,\xi) + N(Y,\xi,\varphi(X)) \, .
\edm
Moreover, in the special case that the Nijenhuis tensor is totally skew-symmetric, we have (see \cite{Friedrich&I}, Lemma 8.3)
\bdm
N(Y, \xi,\varphi(X)) \ = \ - \, d F(\varphi(X), \varphi(Y), \xi) \ = \ 
- \, d F(X, Y, \xi) \, ,
\edm
and the proof of the proposition follows directly.
\end{proof}
\noindent
In dimension five the condition that the Nijenhuis tensor is totally
skew-symmetric is equivalent to its vanishing (see \cite{CG}). 
Therefore we will consider only $5$-dimensional, normal almost contact metric
structures with a Killing vector field $\xi$. The integrability
condition for the existence of a $\nabla$-parallel spinor field in one of the
$1$-dimensional bundles $\Sigma^{\pm}$ was investigated in the paper
\cite{Friedrich&I}. Here we will study the same problem in case that
the $\nabla$-parallel spinor is in $\Sigma^2$. First we suppose that
the normal almost contact structure is regular, i.\,e., the orbit space
$N^4 := M^5/\xi$ is a smooth manifold. In this situation the projection
\bdm
\pi \, : \, M^5 \longrightarrow N^4
\edm
is a principal $S^1$-bundle with curvature form $\hat{\Omega}:=d \eta$. 
Since the Lie derivative $\mathcal{L}_{\xi}\varphi = 0$ of the endomorphism 
$\varphi$ with respect to $\xi$ vanishes (see \cite{Blair}, p.\,49), the 
manifold admits 
\begin{enumerate}
\item a Riemannian metric $\hat{g}$,
\item an integrable complex structure $\hat{\varphi}$ corresponding to 
$\varphi$,
\item a form $\hat{F}(\hat{X},\hat{Y}):= \hat{g}(\hat{X},
\hat{\varphi}(\hat{Y}))$. In general $(N^4, 
\hat{g},\hat{\varphi})$ is not a K\"ahler manifold.
\item a closed $2$-form $\hat{\Omega}$.
\end{enumerate}
The manifold $M^5$ is oriented by the positive frame
$\{e_1, \  e_2 =  - \,\varphi(e_1), \ e_3, \ e_4 =   - \,\varphi(e_3), 
\ e_5 \}$
and we orient the manifold $N^4$ in the same way. Then $\hat{\varphi}$ is a 
section in the {\it positive} twistor bundle
$\mathcal{Z}^+(N^4)$ of the $4$-manifold $N^4$. Suppose that there exists
a $\nabla$-parallel spinor $\psi$,  
\bdm
\nabla^g_X\psi + \frac{1}{4} \cdot (X \haken T) \cdot \psi \ = \ 0, 
\edm
and let us compute its Lie derivative  (see \cite{BG}):
\bdm
\mathcal{L}_{\xi}\psi \ = \ \nabla^g_{\xi}\psi - \frac{1}{4} \cdot d \eta 
\cdot \psi \ = \ - \, \frac{1}{2} \cdot d \eta \cdot \psi\, .
\edm
Since $d \eta (\varphi(X),\varphi(Y)) = d \eta(X,Y)$ (see \cite{Blair}, 
page 51), we can  write $d \eta$ in the form
\bdm
d \eta \ = \ a \cdot e_1 \wedge e_2 + b \cdot (e_1 \wedge e_3 + e_2 \wedge 
e_4) + c \cdot (e_1 \wedge e_4 - e_2 \wedge e_3) + d \cdot e_3 \wedge e_4 \, .
\edm
If the spinor field $\psi$ is a section in the $2$-dimensional bundle
$\Sigma^2$, then $d \eta \cdot \psi = 0$ holds 
if and only if $d \eta$ is proportional
to $F$. The spinor field $\psi$ defines a parallel 
spinor field $\hat{\psi}$ on the manifold $N^4$. Since the 
$\Spin(4)$-representation $\Delta^2_5$ is isomorphic to $\Delta^-_4$, 
the parallel spinor field $\hat{\psi} \in \Gamma(N^4;\Sigma^-)$ is a section 
in the negative spinor bundle and it induces a negative 
complex structure $\hat{\psi} \in \mathcal{Z}^-(N^4)$. 
Denote by $\hat{N}^4$ the manifold $N^4$ equipped with the opposite 
orientation. Then $(\hat{N}^4, \hat{g}, \hat{\psi})$ is a Ricci flat, 
antiself-dual K\"ahler manifold and $\hat{\varphi} \in  
\mathcal{Z}^-(\hat{N}^4)$
becomes a negative, integrable complex structure with closed fundamental form, $d\hat{F} = 0$. Therefore $\hat{\varphi}$ is a parallel complex structure and,
consequently, the negative part $W_-$ of the Weyl tensor vanishes, too. In 
particular, we have proved
\begin{thm}  
Let $(M^5, g, \xi, \eta, \varphi)$ be a $5$-dimensional, compact 
 $\alpha$-Sasakian spin manifold and denote by $\nabla$ the unique 
connection with 
totally skew-symmetric torsion $T= \eta \wedge d \eta$. Moreover, suppose 
that the contact structure is regular. If there exists a
$\nabla$-parallel spinor field $\psi$ such that $F \cdot \psi = 0$, then
the $\alpha$-Sasakian structure is an $S^1$-bundle over a flat torus. 
\end{thm}
\noindent
If the spinor $\psi$ is a section in the $1$-dimensional bundle
$\Sigma^{\pm}$, then $d \eta \cdot \psi = 0$ holds if and only if
$a = - d$. The projected spinor $\hat{\psi}$ is parallel with respect to
a metric connection with torsion, the metric $\hat{g}$ 
is antiself-dual and the $2$-form $\hat{\Omega}$ is antiself-dual, too. 
The space $N^4$ is a hyperhermitian manifold with torsion (or, equivalently,
a hypercomplex manifold). In the compact case there are only three 
possibilities: the flat torus, a $K3$-surface or a Hopf surface 
(see \cite{DI, IP}).
\begin{thm}  
Let $(M^5, g, \xi, \eta, \varphi)$ be a $5$-dimensional, compact, regular 
normal almost contact metric spin manifold with a Killing vector field $\xi$ 
and denote by $\nabla$ the unique connection with 
totally skew-symmetric torsion. If there exists a projectable
$\nabla$-parallel spinor field $\psi$ in the bundle $\Sigma^{\pm}$, then
$M^5$ is a $S^1$-bundle over a flat torus, a $K3$-surface or over a Hopf 
surface. If the projectable $\nabla$-parallel spinor belongs to the 
$2$-dimensional bundle $\Sigma^2$, then $M^5$ is a $S^1$-bundle over a 
flat torus or over a Hopf surface.
\end{thm}
%
\section{The quasi-Sasakian case}\noindent
%
We will integrate the structure equations of a quasi-Sasakian manifold 
admitting a $\nabla$-parallel spinor of type $F \cdot \psi = 0$. Choosing 
the orthonormal 
frame $e_1, \ldots , e_5$ in the $\U(1)$-reduction defined by the 
$\nabla$-parallel spinor, an easy computation yields the proof of the following
\begin{prop}  
Let $(M^5, g, \xi, \eta, \varphi)$ be a normal almost contact metric 
spin manifold with Killing vector field $\xi$. Denote by $\nabla$ the unique 
connection with totally
skew-symmetric torsion T and preserving the structure. Moreover, suppose that
there exists a $\nabla$-parallel spinor field of type $F \cdot \psi = 0$. 
Then, locally, there exists an orthonormal frame as well as a $1$-form $A$
such that the structure equations of the Riemannian manifold are given by
\bdm
de_1 \, = \, A \wedge e_2 +\, e_1 \haken T, \quad  de_2 \, = \, - \, 
A \wedge e_1 +\, e_2 \haken T, \quad de_3 \, = \, A \wedge e_4 +\, e_3 
\haken T, 
\edm
\bdm 
de_4 \, = \, - \, A \wedge e_3 +\, e_4 \haken T ,\quad
de_5 \, = \,  e_5 \haken T \, .
\edm
\end{prop}
\begin{exa}\label{sasaki}
In $\R^5$ we consider the $1$-forms
\bdm
e_1 \, =\, \frac{1}{2} dx_1, \quad e_2 \, =\, \frac{1}{2} dy_1, \quad
e_3 \, =\, \frac{1}{2} dx_2, \quad e_4 \, =\, \frac{1}{2} dy_2,\quad
e_5 \, =\, \eta \, =\, \frac{1}{2} (dz - y_1 \cdot dx_1 - y_2 \cdot dx_2)\, .
\edm
We obtain a Sasakian manifold (see \cite{Blair}) and  it is not 
hard to see that this Sasakian manifold admits
$\nabla$-parallel spinors of type $F \cdot \psi = 0$. The Sasakian structure
arises from left invariant vector fields on a $5$-dimensional Heisenberg group.
\end{exa}
\noindent
First we will prove that Example~\ref{sasaki} is locally the only 
Sasakian manifold admitting a $\nabla$-parallel spinor of type $F \cdot \psi =
0$. Partially, we studied this case already in \cite{Friedrich&I}.
\begin{thm}  
Let $(M^5, g, \xi, \eta, \varphi)$ be a $5$-dimensional Sasakian
spin manifold and denote by $\nabla$ the unique connection with totally
skew-symmetric torsion $T= \eta \wedge d \eta$. If there exists a
$\nabla$-parallel spinor field $\psi$ such that $F \cdot \psi = 0$, then
the Sasakian structure is locally equivalent to the structure described
in Example~$\ref{sasaki}$.
\end{thm}
\begin{proof} 
In case the almost contact metric structure is Sasakian, the
exterior derivative $d \eta$ is proportional to the fundamental form,
\bdm
d \eta \ = \ 2 \cdot F, \quad T \ = \ 2 \cdot (e_1
\wedge e_2 + e_3 \wedge e_4) \wedge e_5 \, .
\edm
The structure equations of the Riemannian manifold can be simplified to:
\bdm
de_1 \ = \ A \wedge e_2 + 2 \cdot e_2 \wedge e_5, \quad
de_2 \ = \ - \, A \wedge e_1 - 2 \cdot e_1 \wedge e_5 \, ,
\edm
\bdm
de_3 \ = \ A \wedge e_4 + 2 \cdot e_4 \wedge e_5, \quad
de_4 \ = \ - \, A \wedge e_3 - 2 \cdot e_3 \wedge e_5 \, ,
\edm
\bdm
de_5 \ = \ 2 \cdot (e_1 \wedge e_2 + e_3 \wedge e_4) \, .
\edm
Differentiating the system we see that the $1$-form $C := A - 2 \cdot e_5$
satisfies the algebraic equations
\bdm
dC \wedge e_1 \ = \ dC \wedge e_2 \ = \ dC \wedge e_3 \ = \ dC \wedge e_4 \ 
= \ 0 \, ,
\edm
i.\,e., $C$ is a closed form, $dC= 0$. Moreover, we obtain
\bdm
\Omega^A \ = \ dA \ = \ 2 \cdot de_5 \ = \ 4 \cdot (e_1 \wedge e_2 + e_3 
\wedge e_4)
\edm
and, in particular, $\Ric^{\nabla} = \diag(-4,-4,-4,-4,0)$. Locally there
exists a function $f$ such that $C$ is the differential of $f$, $C= df$. Using the new frame
\bdm
e_1^* \ := \ \cos(f) \cdot e_1 - \sin(f) \cdot e_2, \quad
e_2^* \ := \ \sin(f) \cdot e_1 + \cos(f) \cdot e_2 \, ,
\edm
\bdm
e_3^* \ := \ \cos(f) \cdot e_3 - \sin(f) \cdot e_4, \quad
e_4^* \ := \ \sin(f) \cdot e_3 + \cos(f) \cdot e_4 \, ,
\edm
we obtain the equations
\bdm
de_1^* \ = \ de_2^* \ = \ de_3^* \ = \ de_4^* \ = \ 0, \quad 
de_5^* \ = \ 2 \cdot (e_1^* \wedge e_2^* + e_3^* \wedge e_4^*) \, .
\qedhere
\edm
\end{proof}
\noindent
We will generalize the argument used for Sasakian manifolds. The aim is
the construction of almost contact metric structures admitting
$\nabla$-parallel spinors of type $F \cdot \psi = 0$. We assume
that the structure under consideration is quasi-Sasakian:
\bdm
N \ = \ 0, \quad dF \ = \ 0, \quad \xi \ \text{is a Killing field} \, .
\edm 
The torsion form is given by $T = \eta \wedge d \eta$ and $d \eta$ has the 
invariance property $d \eta (\varphi(X),\varphi(Y)) = d \eta (X,Y)$ 
(see \cite{Blair}, p.\,$51$). Therefore there exist
functions $a, b, c, d$ such that
\bdm
d \eta \ = \ a \cdot e_1 \wedge e_2 + b \cdot (e_1 \wedge e_3 + e_2 \wedge 
e_4) + c \cdot (e_1 \wedge e_4 - e_2 \wedge e_3) + d \cdot e_3 \wedge e_4 
\edm
holds. The structure equations yield the following system:
\renewcommand{\theequation}{$*$}  
\bea
\notag
de_1 & = & A \wedge e_2 + a \cdot e_2 \wedge e_5 + b \cdot e_3 \wedge e_5 +
c \cdot e_4 \wedge e_5 \, ,\\
\notag
de_2 & = & - \, A \wedge e_1 - a \cdot e_1 \wedge e_5 + b \cdot e_4 \wedge e_5
- c \cdot e_3 \wedge e_5 \,,\\
de_3 & = & A \wedge e_4 - b \cdot e_1 \wedge e_5 + c \cdot e_2 \wedge e_5 +
d \cdot e_4 \wedge e_5 \,,\\
\notag
de_4 & = & - \, A \wedge e_3 - c \cdot e_1 \wedge e_5 - b \cdot e_2 \wedge e_5
- d \cdot e_3 \wedge e_5 \,,
\eea
\bdm
de_5  =  a \cdot e_1 \wedge e_2 + b \cdot (e_1 \wedge e_3 + e_2 \wedge 
e_4) + c \cdot (e_1 \wedge e_4 - e_2 \wedge e_3) + d \cdot e_3 \wedge e_4 \, .
\edm
We suppose now that $a, b, c, d$ are constant. For example, this assumption is
satisfied if the quasi-Sasakian structure is locally homogeneous. We study 
the integrability condition of the system ($*$). A straightforward, but 
lengthy computation yields that there is only one integrability condition. 
\begin{lem} 
Let $a, b, c, d$ be constant and $A$ be a $1$-form. The system 
$(*\,)$ is
integrable if and only if the equation
\bdm
\Omega^A \ = \ dA \ = \ (ad - b^2 - c^2) \cdot F
\edm
holds.
\end{lem}
\noindent
Conversely, fix real parameters $a, b, c, d$ and suppose that $b^2 + c^2 -
ad \neq 0$. Then there exist (see \cite{Fla}, chapter 7.6) locally
$1$-forms $e_1, \ldots , e_5, A$ solving the system ($*$) and the equation
\bdm
 dA \ = \ (ad - b^2 - c^2) \cdot F \ .
\edm
We define a Riemannian metric as well as an almost contact metric structure
by the condition that $e_1,\,\ldots, e_5 = \eta$ is an orthonormal frame. 
We denote the corresponding almost contact metric manifold  by $M^5(a,b,c,d)$.
Its connection forms $\omega_{ij}$ of the Levi-Civita connection are:
\bdm\ba{cccccc}
2 \cdot \omega_{15}&\!\!\!=\!\!\!&(+a \cdot e_2 + b \cdot e_3 + c \cdot e_4),&
2 \cdot \omega_{25}&\!\!\!=\!\!\!&(-a \cdot e_1- c \cdot e_3+ b \cdot e_4),\\
2 \cdot \omega_{35}&\!\!\!=\!\!\!&(- b \cdot e_1 + c \cdot e_2 + d \cdot e_4),&
2 \cdot \omega_{45}&\!\!\!=\!\!\!&(- c \cdot e_1 - b \cdot e_2 - d \cdot e_3),\ea
\edm
\vspace{-4mm}
\bdm\ba{cccccrccc}
2 \cdot \omega_{12}&\!\!=\!\!&2 \cdot A - a \cdot e_5,& 
2 \cdot \omega_{13}&\!\!=\!\!&- \, b \cdot e_5,&
2 \cdot \omega_{14}&\!\!=\!\!&- \, c \cdot e_5,\\
2 \cdot \omega_{34}&\!\!=\!\!& 2 \cdot A  - \, d \cdot e_5,&
2 \cdot \omega_{23}&\!\!=\!\!&c \cdot e_5,& 
2 \cdot \omega_{24}&\!\!=\!\!&- \, b \cdot e_5\,.
\ea\edm
We verify now by a direct computation that $\xi$ is a Killing vector field,
the Nijenhuis tensor vanishes and the fundamental form $F$ is
closed, $dF = 0$. Moreover, the torsion form $T = \eta \wedge d \eta$ is
coclosed, $\delta(T) = 0$. Let us summarize the result.
\begin{thm}
The almost contact metric manifold $M^5(a,b,c,d)$ has the following
properties:
\begin{enumerate}
\item $\xi$ is a Killing vector field.
\item The Nijenhuis tensor vanishes, $N = 0$.
\item The fundamental form is closed, $dF = 0$.
\item $d \eta \wedge d \eta = ( a d - b^2 - c^2) \cdot F \wedge F$.
\item The torsion form $T = \eta \wedge d \eta$ is coclosed, $\delta(T) = 0$.
\item There exist two $\nabla$-parallel spinor fields with $F \cdot \psi = 0$.
\item The Ricci tensor of the connection $\nabla$ is 
\bdm
\Ric^{\nabla} \ = \ (b^2 + c^2 - ad) \cdot \diag(1,\, 1,\, 1,\, 1,\, 0)\,.
\edm
\end{enumerate}
Conversely, if $M^5$ is a locally homogeneous, quasi-Sasakian 
spin manifold with a $\nabla$-parallel spinor of type $F \cdot \psi = 0$, then 
$M^5$ is locally equivalent to one of the manifolds $M^5(a,b,c,d)$. 
\end{thm}
\noindent
For any metric connection with totally skew-symmetric torsion the action of
the $4$-form
\bdm
2 \cdot \sigma^T \ : = \ \sum_{i=1}^{5} (e_i \haken T) \wedge (e_i \haken T)
\edm
on spinors plays an important role (see \cite{Friedrich&I}). For the 
quasi-Sasakian structures  $M^5(a,b,c,d)$ we immediately obtain 
the formula
\bdm
dT + 2 \cdot \sigma^T + \Scal^{\nabla} \ = \ 4 \cdot (ad - b^2 - c^2) \cdot 
(e_1 \wedge e_2 \wedge e_3 \wedge e_4 \, - \, 1) \, .
\edm
Moreover, on $\Sigma^2$ the relation $e_1 \cdot e_2 = - \, e_3 \cdot 
e_4$ holds and then we obtain
\bdm
e_1 \cdot e_2 \cdot e_3 \cdot e_4 \, - \, 1 \ = \ 0 \, .
\edm 
Finally we apply a general integral formula (see \cite{Friedrich&I}).
\begin{thm} 
Let $(M^5,g, \xi, \eta, \varphi)$ be a compact, locally homogeneous 
quasi-Sasakian spin manifold. 
Denote by $D^{\nabla}$ the Dirac operator defined by
its unique connection $\nabla$ with totally skew-symmetric torsion. Then any
$\nabla$-harmonic spinor of type $F \cdot \psi = 0$ is $\nabla$-parallel,
\bdm
\{D^{\nabla} \psi = 0, \ F \cdot \psi = 0 \} \ = \ \{ \nabla \psi = 0, \ 
F \cdot \psi = 0 \} \, .
\edm
\end{thm} 
\noindent
Remark that $dT + 2 \cdot \sigma^T + \Scal^{\nabla}$ acts on the bundles
$\Sigma^{\pm}$ by a multiplication by $\pm 8 \cdot (ad - b^2 - c^2)$. Hence, 
in one of these $1$-dimensional subbundles of the spinor bundle there are no
$\nabla$-harmonic spinors at all.\\

\noindent
We will interpret the quasi-Sasakian structures $M^5(a,b,c,d)$. The system
\bea[*]
dA &= & (ad - b^2 - c^2) \cdot (e_1 \wedge e_2 + e_3 \wedge e_4) \, , \\
de_1 & = & A \wedge e_2 + a \cdot e_2 \wedge e_5 + b \cdot e_3 \wedge e_5 +
c \cdot e_4 \wedge e_5 \, ,\\
de_2 & = & - \, A \wedge e_1 - a \cdot e_1 \wedge e_5 + b \cdot e_4 \wedge e_5
- c \cdot e_3 \wedge e_5 \,,\\
de_3 & = & A \wedge e_4 - b \cdot e_1 \wedge e_5 + c \cdot e_2 \wedge e_5 +
d \cdot e_4 \wedge e_5 \,,\\
de_4 & = & - \, A \wedge e_3 - c \cdot e_1 \wedge e_5 - b \cdot e_2 \wedge e_5
- d \cdot e_3 \wedge e_5 \,,
\eea[*]
\vspace{-5mm}
\bdm
de_5  \ = \  a \cdot e_1 \wedge e_2 + b \cdot (e_1 \wedge e_3 + e_2 \wedge 
e_4) + c \cdot (e_1 \wedge e_4 - e_2 \wedge e_3) + d \cdot e_3 \wedge e_4
\edm
has a solution in $\R^6$ (see \cite{Fla}, Lie's third Theorem). Moreover, 
there exists a $6$-dimensional Lie group $G(a,b,c,d)$ such that the 
forms $e_1, \ldots , e_5, A$ constitute a basis of the left invariant forms. 
Any submanifold $M^5 \subset G(a,b,c,d)$ such that the restricted
forms $e_1, \ldots , e_5$ are independent is a model for the 
quasi-Sasakian structure $M^5(a,b,c,d)$. In the Sasakian case 
$(a = d = 2$, $c = d = 0)$ the Lie group $G(1,0,0,1)$ is the product of the 
$5$-dimensional Heisenberg group by $\R^1$. Let us discuss the case of
$ a = c = d = 0$, $b =1$. Introducing the $1$-forms
\bdm
\sqrt{2} \cdot e^*_1 \ := \ e_1 + e_4, \quad
\sqrt{2} \cdot e^*_2 \ := \ e_2 - e_3, \quad
\sqrt{2} \cdot e^*_3 \ := \ e_1 - e_4, \quad
\sqrt{2} \cdot e^*_4 \ := \ e_2 + e_3, \quad
\edm
\bdm
C^*_1 \ := \ A + e_5, \quad C^*_2 \ := \ A - e_5\,,
\edm
this system is equivalent to
\bdm
de^*_1 \ = \ C^*_1 \wedge e^*_2, \quad
de^*_2 \ = \ - \, C^*_1 \wedge e^*_1, \quad
de^*_3 \ = \ C^*_2 \wedge e^*_4, \quad
de^*_4 \ = \ - \, C^*_2 \wedge e^*_3,
\edm
\bdm
d C^*_1 \ = \ - \, 2 \cdot e^*_1 \wedge e^*_2, \quad 
d C^*_2 \ = \ - \, 2 \cdot e^*_3 \wedge e^*_4 \, .
\edm
Consequently, in case of $a = c = d = 0$, $b =1$, the group $G$ is isomorphic
to the group $S^3 \times S^3$ and a generic submanifold of codimension
one admits a quasi-Sasakian structure of that type. An interesting example
is the $5$-dimensional Stiefel manifold $V_{4,2}$ defined as 
the set of orthogonal pairs $(x,y)$ of vectors in $\R^4$ of length one,
\bdm
V_{4,2} \ = \ \big\{(x,y) \in \R^4 \, : \ ||x|| \, = \, ||y|| \, = \, 1, \ 
g(x,y) \, = \, 0 \big\} \ \subset \ S^3 \times S^3.
\edm
The Stiefel manifold is a naturally reductive space, it admits an Einstein 
metric constructed first by J. Jensen (see \cite{Jen}) as well as real 
Killing spinors (see \cite{Fri1}). The Killing spinor induces a Sasakian structure on $V_{4,2}$, but the spinor does not belong to the subbundle 
$\Sigma^2$ with respect to this Sasakian structure. However, there is a 
second quasi-Sasakian 
structure realizing our parameters $a = c = d = 0$ and $ b \neq 0$ and it corresponds to the embedding of $V_{4,2}$ into $S^3 \times S^3$. Finally one can
construct a third, nonquasi-Sasakian structure on $V_{4,2}$ with
\bdm
N \ = \ 0, \quad dF \ \neq \ 0, \quad d^{\varphi}F \ = \ 0 \, .
\edm 
Its torsion is thus again given by $T = \eta \wedge d \eta$, it admits a
$\nabla$-parallel spinor, but this spinor does not belong to the bundle
$\Sigma^2$. In a forthcoming paper (see \cite{Agri}) this example and more general connections on it will be discussed from the point of view of naturally 
reductive spaces. In particular, it turns out that for all three structures
the corresponding connection coincides with the canonical connection of the reductive space.\\

\noindent
In the degenerate case, i.e., $a = 1$, $b = c = d = 0$, we have to solve 
\bdm
de_1 \ = \ A \wedge e_2 + e_2 \wedge e_5, \quad
de_2 \ = \ - \, A \wedge e_1 - e_1 \wedge e_5, \quad
de_3 \ = \ A \wedge e_4, \quad
de_4 \ = \ - \, A \wedge e_3, 
\edm
\vspace{-5mm}
\bdm
de_5 \ = \ e_1 \wedge e_2, \quad dA \ = \ 0 \, . 
\edm 
Locally the $1$-form $A$ is the differential of some function, $A = df$. We 
introduce again a new frame $e_1^*, \ldots , e_5^*$ by the same formulas as 
in the proof of Theorem 4.1. The equations now read
\bdm
de^*_1 \ = \ e^*_2 \wedge e^*_5, \quad
de^*_2 \ = \ e^*_5 \wedge e^*_1, \quad
de^*_5 \ = \ e^*_1 \wedge e^*_2, \quad
de^*_3 \ = \ 0, \quad de^*_4 \ = \ 0 
\edm
and, consequently, the degenerate quasi-Sasakian structure $M^5(1,0,0,0)$ is 
realized on the Lie group $S^3 \times \R^2$.
%
\section{The dilation function and conformal transformations}\noindent
%
Consider an almost contact metric structure
$(M^{2k+1}, g, \xi , \eta , \varphi)$ with skew-symmetric Nijenhuis
tensor $N$ and Killing vector field $\xi$. 
The torsion of its unique connection $\nabla$ is given by
\bdm
T \ = \ \eta \wedge d \eta + d^{\varphi}F + N - \eta \wedge (\xi \haken N) \, .
\edm 
Let us introduce the \emph{Lee form} $\theta$ defined by the formula
\bdm
\theta(X) \ := \ - \, \frac{1}{2} \sum_{i=1}^{2k+1} T(\varphi(X),e_i, \varphi(e_i)) \, .
\edm
Remark that all vectors in the tuple $(\varphi(X), e_i, \varphi(e_i))$ are 
orthogonal to $\xi$ except in case  $e_i = \xi$. Consequently, we have
\bdm
\eta \wedge d \eta (\varphi(X), e_i, \varphi(e_i)) \ = \ 0 \, , \quad
\eta \wedge (\xi \haken N) (\varphi(X), e_i, \varphi(e_i)) \ = 0 \, .
\edm
Moreover, since $N$ is skew-symmetric, we obtain the formula
\bdm
N\big(\varphi(X), e_i, \varphi(e_i)\big) \ = \ g\big(N(e_i,\varphi(e_i)), 
\varphi(X)\big) \ = 
\ d \eta\big(e_i,\varphi(e_i)\big)\cdot g\big(\xi, \varphi(X)\big) \ = \ 0 \, .
\edm
To summarize, the Lee form $\theta$ depends only on the exterior derivative of 
the fundamental form,
\bdm
\theta(X) \ := \ - \, \frac{1}{2} \sum_{i=1}^{2k+1} dF(X,e_i, \varphi(e_i)) 
\edm
We define a new almost contact metric structure by the formulas
\bdm
\varphi' \ := \ \varphi, \quad \xi' \ := \ \xi \quad \eta' \ := \ \eta,
\quad g' \ := \ e^{2f} \cdot g + (1 - e^{2f}) \cdot \eta \otimes \eta \, .
\edm
If the function $f$ is constant along the integral curve of $\xi$, $df(\xi) 
= 0$,  we call the transformation \emph{special conformal}. Moreover,
computing the Lie derivative $\mathcal{L}_{\xi}g'$ of the new metric $g'$ 
with respect to the vector field, we obtain the following
\begin{prop}
The class of all almost contact metric manifolds with skew-symmetric 
Nijenhuis tensor and Killing vector field $\xi$ is invariant under special
conformal transformations.
\end{prop}
\noindent
In particular, if $(M^{2k+1}, g, \xi, \eta, \varphi)$ admits
a connection with totally skew-symmetric torsion and preserving the structure,
then the transformed manifold admits again such a  connection. Let us 
compute the new torsion form $T'$ as well as the new Lee form $\theta'$. 
Since the fundamental form and the Nijenhuis tensor transform as
\bdm
F' \ = \ e^{2f} \cdot F \, , \quad  N' \ = \ N
\edm
and the new torsion form is given by $T'=\eta'\wedge d\eta'+d^{\varphi'}F' + 
N' - \eta' \wedge (\xi' \haken N')$, 
we have
\begin{prop} 
Let $(M^{2k+1}, g, \xi, \eta, \varphi)$ be an almost contact 
metric manifold with skew-symmetric Nijenhuis tensor and Killing vector 
field $\xi$. The torsion form $T'$ and the Lee form $\theta'$ of a 
structure obtained by a special conformal transformation are given 
by the formulas:
\bdm
T' \ = \ T + (e^{2f} - 1) \cdot d^{\varphi} F - 2 \cdot 
e^{2f} \cdot (df \circ \varphi) \wedge F, \quad \theta' \ = \ \theta + 
2 \cdot df \, .
\edm
\end{prop}
\noindent
We compute the relation between the spinorial covariant derivatives 
$\nabla^{g}$ , $\nabla^{g'}$ corresponding to the Levi-Civita connections
of $g$ and $g'$ as well as between $\nabla$ and $\nabla'$. Let us summarize 
the result:
\begin{prop} Let $X$ be a vector field orthogonal to $\xi$ and $\psi$ be an arbitrary spinor field. Then the following formulas hold:
\bdm
\nabla^{g'}_X \psi \ = \ \nabla^{g}_X \psi + \frac{1}{4} \cdot (df \cdot X - 
X \cdot df) \cdot \psi + \frac{1 - e^{-f}}{4} \cdot (X \haken d \eta) \cdot 
\xi \cdot \psi \, ,
\edm
\bdm
\nabla^{g'}_{\xi} \psi \ = \ \nabla^{g}_{\xi}\psi + \frac{e^{-2f} - 1}{4} 
\cdot d \eta \cdot \psi \, .
\edm
If, in addition, the Nijenhuis tensor $N = 0$ vanishes, one has 
\bdm
\nabla'_X \psi \ = \ \nabla_X \psi + \frac{1}{4} \cdot (df \cdot X - 
X \cdot df) \cdot \psi - \frac{1}{2} \cdot X \haken (df \circ \varphi) 
\wedge F \cdot \psi \, ,
\edm
\bdm
\nabla'_{\xi} \psi \ = \ \nabla_{\xi}\psi + \frac{e^{-2f} - 1}{2} \cdot d \eta \cdot \psi \,.
\edm
\end{prop}
\begin{proof} The first two formulas are direct consequences of the
transformation of the Levi-Civita connection. For example, let us compute 
the third formula. First remark that the new torsion form $T'$ can be 
expressed as
\bdm
T' \ = \ e^{2f} \cdot T + (1 - e^{2f}) \cdot \eta \wedge d \eta - 
2 \cdot e^{2f} \cdot (df \circ \varphi) \wedge F \, .
\edm
In a fixed $g'$-orthonormal frame $e'_1 := e^{-f} \cdot e_1, \ldots , 
e'_{2k} := e^{-f} \cdot e_{2k}, e'_{2k+1} : = e_{2k+1} = \xi$ we obtain
\bdm
\nabla'_X \psi \ = \ \nabla^{g'}_X \psi + \frac{1}{8} \sum^{2k}_{i,j=1}
T'(X,e'_i,e'_j) \cdot e_i \cdot e_j \cdot \psi +
\frac{1}{4} \sum_{i=1}^{2k}
T'(X,e'_i,\xi) \cdot e_i \cdot \xi \cdot \psi  
\edm
\bdm
= \ \nabla^{g'}_X \psi + \frac{1}{8} \sum^{2k}_{i,j=1}
T(X,e_i,e_j) \cdot e_i \cdot e_j \cdot \psi - 2 \cdot 
\big(X \haken (df \circ \varphi) \wedge F\big) \cdot \psi +
\frac{e^{-f}}{4} \cdot (X \haken d \eta) \cdot \xi \cdot \psi \, .
\edm
We used that $\eta \wedge d \eta (X,e_i, e_j)= 0$ vanishes for $X$ 
orthogonal to $\xi$ and $i,j \leq 2k$ as well as the formula
\bdm
T'(X, e'_i, \xi) \ = \ (\xi \haken T')(X, e'_i) \ = \ d \eta (X, e'_i) \ = \ 
e^{-f} \cdot d \eta(X, e_i) \, .
\edm
\noindent
We apply now the formula for the relation between the spinorial connections $\nabla^{g'}$ and $\nabla^g$ and combine the latter formula with the definition 
of the connection $\nabla$. Then the term involving 
$(X \haken d \eta) \cdot \xi \cdot \psi$ cancels and finally we get 
the required formula.
\end{proof}
\noindent
In dimension five one verifies directly that the following endomorphisms
acting on the different parts of the spinor coincide:
\bdm
 2 \cdot X \haken (df \circ \varphi) \wedge F \ = \ df \cdot X - X \cdot df 
\quad \text{on the bundle} \quad \Sigma^{\pm} \, ,
\edm
\bdm
 2 \cdot X \haken (df \circ \varphi) \wedge F \ = \ - df \cdot X + X \cdot df 
\quad \text{on the bundle} \quad \Sigma^{2} \, .
\edm
Consequently, the formulas for the transformation of the spinorial covariant
derivative in directions $X$ orthogonal to $\xi$ can be simplified,
\bdm
\nabla'_X \psi \ = \ \nabla_X \psi \quad \text{for} \  
\psi \in \Gamma(\Sigma^{\pm}), \quad 
\nabla'_X \psi \ = \ \nabla_X \psi + \frac{1}{2}(df \cdot X - X \cdot df) \cdot \psi \quad \text{for} \ \psi \in \Gamma(\Sigma^{2}).
\edm
We see that a $\nabla$-parallel spinor in the bundle $\Sigma^2$ never transforms into a $\nabla'$-parallel spinor excepted if $f$ is constant. In the bundle 
$\Sigma^{\pm}$ the situation is different, it may occur.
\begin{thm}
Let $(M^5, g, \xi, \eta, \varphi, \nabla)$ and  $(M^5, g', \xi, \eta, \varphi, 
\nabla')$ be two $5$-dimensional almost contact metric normal manifolds with
Killing vector field $\xi$ which are connected by a special conformal 
transformation. If $d \eta$ is antiself-dual then a spinor is $\nabla$-parallel
if and only if it is $\nabla'$-parallel.
\end{thm}
\noindent
In dimension five the Lee form determines completely the exterior derivative 
$dF$ of the fundamental form, provided the vector field $\xi$ is Killing:
\bdm
d F \ = \ \theta \wedge F, \quad d^{\varphi} F \ = \ - \, 
(\theta \circ \varphi) \wedge F, \quad T \ = \ \eta \wedge d \eta - 
(\theta \circ \varphi) \wedge F \, .
\edm
Consider a structure with closed Lee form $\theta$ and solve locally the
equation $2 \cdot df + \theta = 0$. Since $\theta(\xi) = 0$, the function $f$
does not depend on $\xi$ and we obtain
\bdm
d(F') \ = \ d(e^{2f} \cdot F) \ = \ e^{2f} \cdot (2 \cdot df \wedge F + dF) \ 
= \ e^{2f} \cdot (2 \cdot df + \theta) \wedge F \ = \ 0 \, .
\edm
\begin{prop} 
Let $(M^5, g, \xi, \eta, \varphi)$ be a $5$-dimensional, normal 
almost contact metric manifold with closed Lee form and Killing vector field $\xi$. Then the space is locally special conformal
to a quasi-Sasakian manifold.
\end{prop}
\noindent
We investigate the algebraic equation $(2 \cdot d \Phi - T) \cdot \psi \, = \, 0$ for the differential of a dilation function $\Phi$ on a $5$-dimensional 
manifold. If the almost contact structure is regular this second Killing 
equation means that the spinor is projectable. Again, there are two cases 
depending on the 
algebraic type of the spinor. The corresponding relations between the $1$-form $d \Phi$ and the $3$-form $T$ were computed in the paper \cite{Friedrich&I} (Lemma 7.2 and lemma 7.5). In order to formulate these relations in an invariant 
way, let us introduce the Hodge operator
$*_{4}$ acting in the $4$-dimensional orthogonal complement of the vector 
field $\xi$.
\begin{prop}
The equation $(2 \cdot d \Phi - T) \cdot \psi \ = \ 0$
\begin{enumerate}
\item admits a solution $0 \neq \psi \in \Sigma^{\pm}$ if and only if
 $\ \ 2 \cdot d \Phi \ = \ - \, \theta, \ *_{4} d \eta \ = \ - \, d \eta \, .$ 
\item admits a solution $0 \neq \psi \in \Sigma^{2}$ if and only if
 $\ \ 2 \cdot d \Phi \ = \ \theta, \ *_{4} d \eta \ = \ d \eta \, .$
\end {enumerate}
\end{prop}
\begin{thm} Let $(M^5, g, \xi, \eta, \varphi)$ be a normal almost contact
metric manifold with Killing vector field $\xi$. If there exist a function
$\Phi$ and an arbitrary, nontrivial spinor $\psi$ such that
\bdm
(2 \cdot d \Phi - T) \cdot \psi \ = \ 0 \, ,
\edm
then the Lee form is closed and $M^5$ is locally special conformal to 
a quasi-Sasakian manifold.
\end{thm}
\noindent
We discuss the case that $\psi$ is a section in the $2$-dimensional 
bundle $\Sigma^2$ in more detail. Since the Nijenhuis tensor vanishes, 
the $2$-form $d \eta$ has the invariance property (see \cite{Blair}, p.\,$51$)
\bdm
d \eta (\varphi(X),\varphi(Y)) \ = \ d \eta (X,Y)\,.
\edm
Together with $*_{4} d \eta = d \eta\ $ the latter relation yields that 
$d \eta$ is proportional to $F$,
\bdm
d \eta \ = \ a \cdot F, \quad 2 \cdot d \Phi \ = \ \theta \, .
\edm
The function $a$ is related to $\Phi$. Indeed, we have
\bdm
0 \ = \ d d \eta \ = \ d a \wedge F + a \cdot d F \ = \ 
(d a + a \cdot \theta) \wedge F \ = \ 
(d a + 2 \cdot a \cdot d \Phi) \wedge F \, ,
\edm
and with respect to $F = e_1 \wedge e_2 + e_3 \wedge e_4$ we conclude that
\bdm
d a + 2 \cdot a \cdot d \Phi \ = \ 0, \quad a \ = \ \frac{1}{C} \cdot 
e^{-2 \Phi} \, .
\edm
We apply a special conformal transformation using the function $f = - \, \Phi$.
The new metric becomes
\bdm
g' \ = \ C \cdot a \cdot g + (1 - C \cdot a) \cdot \eta \otimes \eta, \quad
F' \ = \ C \cdot a \cdot F \ = \ C \cdot d \eta' \, 
\edm
i.\,e., the manifold $(M^5, g', \xi', \eta', \varphi')$ is $\alpha$-Sasakian.
\begin{thm} 
Let $(M^5, g , \xi , \eta , \varphi)$ be a normal almost contact metric 
spin manifold with a Killing vector $\xi$ and let $\Phi$ be a smooth function. 
If there exists an arbitrary, nontrivial spinor field 
$\psi \in \Gamma(\Sigma^2)$ such that
\bdm
(2 \cdot d \Phi - T) \cdot \psi \ = \ 0 \, ,
\edm
then $M^5$ is specially conformal to a Sasakian manifold  $(M^5, g' , \xi , 
\eta , \varphi)$. The metric $g'$ as well as the torsion $T'$ of the connection
$\nabla'$ are given by the formulas
\bdm
g' \ = \ e^{-2\Phi} \cdot g + (1 - e^{-2\Phi}) \cdot \eta \otimes \eta, \quad
T' \ = \ T + 2 \cdot (d \Phi \circ \varphi) \wedge F \, .
\edm
The spinor field $\psi$ is a solution of the algebraic equation $T' \cdot 
\psi = 0$ and for any spinor field in $\Sigma^2$ we have
\bdm
\nabla'_X \psi \ = \ \nabla_X \psi + \frac{1}{2} \cdot (X \cdot d \Phi - 
d \Phi \cdot X) \cdot \psi, \quad \nabla'_{\xi} \psi \ = \ \nabla_{\xi} 
\psi \, .
\edm 
\end{thm}
\noindent
We did not assume that the spinor field $\psi$ is $\nabla$-parallel. But now let us consider a $\nabla$-parallel spinor field $\psi$ in $\Sigma^2$ on 
a quasi-Sasakian spin manifold,
\bdm
\nabla \psi \ = \ 0, \quad (2 \cdot d \Phi - T) \cdot \psi \ = \ 0, 
\quad dF \ = \ 0 \, .
\edm
The relation $d \eta = 
a\cdot F$ yields that the function $a$ is constant and, consequently,
the dilation function $\Phi$ is constant, too.
\begin{thm} 
Let $(M^5, g , \xi , \eta , \varphi)$ be a quasi-Sasakian 
spin manifold with a Killing vector $\xi$ and let $\Phi$ be a smooth function. 
If there exists a nontrivial $\nabla$-parallel spinor field 
$\psi \in \Gamma(\Sigma^2)$ such that
\bdm
(2 \cdot d \Phi - T) \cdot \psi \ = \ 0 \, , \quad \nabla \psi \ = \ 0 \, ,
\edm
then the dilation function $\Phi$ is constant and $M^5$ is special 
conformal homothetic to  a Sasakian manifold of Example $\ref{sasaki}$.
\end{thm}
\begin{exa}
Consider the manifolds $M^5(a,b,c,d)$ . Using the 
formula for the $2$-form $d \eta$ as well as the matrix representation of the 
$5$-dimensional Clifford algebra we compute that the endomorphism 
$T = \eta \wedge d \eta$ acts on $\Sigma^2$ as a symmetric endomorphism with 
two eigenvalues ${\pm} \sqrt{(a-d)^2 + 4 b^2 + 4 c^2}$. In particular, a 
$\nabla$-parallel spinor satisfying the algebraic equation
$T \cdot \psi = 0$ exists only in case that $a=d$ and $b=0, c=0$.
\end{exa}
%
%
\section{Harmonic $1$-forms in the presence of $\nabla$-parallel spinors}
\noindent
%
A classical Theorem of Tachibana states that on a compact $K$-contact 
manifold $(M^{2k+1}, g, \xi, \eta, \varphi)$ any harmonic $1$-form $\omega$ is 
orthogonal to $\eta$, $\ \xi \haken \omega = 0$ (see \cite{Blair}, page 69). 
The $K$-contact condition means by definition that $\xi$ is a Killing vector 
field and $d \eta$ is proportional to the fundamental form $F$, $d \eta = 2 \cdot F$. In particular, in Tachibana's Theorem there is no assumption about the
Nijenhuis tensor $N$. We will prove a similar result, but in a quite different
situation. 
\begin{thm}
Let $(M^{2k+1}, g, \xi, \eta, \varphi)$ be a compact, almost contact metric
manifold with totally skew-symmetric Nijenhuis tensor $N$ and Killing vector
field $\xi$. Denote by $T$ the torsion form of its contact unique connection 
$\nabla$ and suppose that
\bdm
\xi \haken \delta \big(T - d^{\varphi}F \big) \ = \ 0
\edm
holds. Then any harmonic $1$-form is orthogonal to $\eta$.
\end{thm}
\begin{proof} Basically, we follow the proof of Tachibana's Theorem. A 
harmonic $1$-form $\omega$ on a compact Riemannian
manifold is invariant under any Killing vector field, $\mathcal{L}_{\xi} \omega = 0$. Since
\bdm
0 \ = \ \mathcal{L}_{\xi} \omega \ = \ \xi \haken d \omega + d(\xi \haken \omega) \ = 
\ d(\xi \haken \omega)
\edm
the function $\omega(\xi)$ is constant. We introduce the $1$-form $\beta := 
\omega - \omega(\xi) \cdot \eta$ and compute
\bdm
\Delta(\beta) \ = \ - \, \omega(\xi) \cdot \Delta(\eta) \ = \ - \, \omega(\xi) 
\cdot \big( \delta d (\eta) + d \delta (\eta)\big) \ = \ - \, \omega(\xi) 
\cdot \delta d (\eta) \, .
\edm
On the other hand, since $\xi$ is a Killing vector field, we have
\bdm
\delta d \eta(X) - \frac{1}{2} \cdot \eta(X) \cdot |d \eta|^2 \ = \ 
\delta d \eta(X) - 2 \cdot \eta(X) \cdot |\nabla^g \eta|^2 \ = \ 
\delta(\eta \wedge d \eta)(X, \xi) \ = \ \delta (T - d^{\varphi} F)(X, \xi) 
\edm
and the assumption yields the equality
\bdm
2 \cdot \delta d (\eta) \ = \ |d \eta|^2 \cdot \eta \, .
\edm
Using that $\beta$ and $\eta$ are orthogonal, $(\beta, \eta) = 0$, we obtain
\bdm
2 \cdot \int_{M^{2k+1}} \big(\Delta(\beta), \beta \big) \ = \ - \, 
\int_{M^{2k+1}} \omega(\xi) \cdot |d \eta|^2 \cdot (\beta, \eta) \ = \ 0
\edm
and, consequently, $\beta$  is a harmonic $1$-form. Since $d \eta$ is a nonvanishing $2$-form we conclude that $\omega(\xi) = 0$.
\end{proof} 
\begin{NB}We proved in fact that any harmonic $1$-form satisfying the
algebraic condition $\delta(T - d^{\varphi}F)(\omega, \xi) = 0$ is orthogonal
to $\eta$.
\end{NB}
\begin{exa} The latter Theorem applies for quasi-Sasakian manifolds 
($N = 0$, $ d F = 0$ and $\xi$ Killing) with divergence free torsion form $T$. Locally homogeneous examples of this type have been discussed in section $4$.
\end{exa}
\begin{exa}\label{exa-qs}
Consider a quasi-Sasakian spin manifold. The spinor bundle $\Sigma$ of 
$M^{2k+1}$ splits into several parts 
and precisely two subbundle $\Sigma^{\pm}$ are $1$-dimensional. If there exists
a $\nabla$-parallel spinor in one of these bundles, then the condition of the
Theorem is satisfied. Indeed, using Proposition $9.1$ of \cite{Friedrich&I}, 
the integrability condition of the parallel spinor yields the equation
\bdm
\delta T(X, \xi) \ = \ - \, \mbox{Ric}^{\nabla}(\xi, X) \ = \ 0 \, .
\edm
\end{exa}
\noindent
Let us again discuss the $5$-dimensional quasi-Sasaki case with a $\nabla$-parallel spinor in more details. If the spinor is a section in one of the bundles
$\Sigma^{\pm}$, then any harmonic $1$-form is orthogonal to $\eta$ (see 
Example \ref{exa-qs}). Suppose next that the $\nabla$-parallel spinor is in $\Sigma^2$. The Ricci tensor of the connection $\nabla$ is given by the curvature form $\Omega^A$ of the connection $A$ (see Theorem \ref{thm-curv}):
\bdm
\mbox{Ric}^{\nabla}(Y, Z) \ = \ \Omega^A(Y , \varphi(Z)) \, .
\edm
From the formula (see \cite{Friedrich&I}, section 2)
\bdm
\delta(T)(X, \xi) \ = \ \mbox{Ric}^{\nabla}(\xi, X) - \mbox{Ric}^{\nabla}(X, 
\xi) \ = \ \Omega^A(\xi, \varphi(X)) - \Omega^A(X, \varphi(\xi)) \ = \ 
\Omega^A(\xi, \varphi(X))
\edm
we see that $\xi \haken \delta(T)$ coincides with $(\xi \haken \Omega^A) \circ \varphi$. We compute $\Omega^A \wedge F$. In fact, we have the general formula 
(see \cite{Friedrich&I})
\bdm
\Omega^A \wedge F \ = \ \hat{\sigma}(R^{\nabla}) \ = \ \frac{3}{2} \cdot d T 
- \sigma^T, \quad \sigma^T := \frac{1}{2} \sum_{i=1}^5 (e_i \haken T) \wedge 
(e_i \haken T)
\edm
and the structure equations for a quasi-Sasakian manifold with parallel
spinor field in $\Sigma^2$ yield 
\bdm
\xi \haken d T \ = \ \xi \haken \sigma^T \ = \ 0 \, .
\edm
Finally, we conclude that $\xi \haken (\Omega^A \wedge F)= 0$ and the form
$\xi \haken \delta(T)$ vanishes.
\begin{thm}
Let $(M^{5}, g, \xi, \eta, \varphi)$ be a compact, $5$-dimensional 
quasi-Sasakian manifold with a $\nabla$-parallel spinor. Then any harmonic $1$-form is orthogonal to 
$\eta$.
\end{thm}
    
\end{document}